\documentclass[11pt]{article}
\newtheorem{theorem}{Theorem}[section]
\newtheorem{corollary}{Corollary}[section]

\usepackage{amssymb}
\usepackage{amsmath}
\usepackage{epsfig}

\begin{document}
\begin{center}
{\LARGE\bf  % 138cor/ar2.tex \\ \today \\
The distribution of the maximum}\\[1ex]
{\LARGE\bf of a second order autoregressive process:}\\[1ex]
{\LARGE\bf the continuous case}\\[1ex]
by\\[1ex]
Christopher S. Withers\\
Applied Mathematics Group\\
Industrial Research Limited\\
Lower Hutt, NEW ZEALAND\\[2ex]
Saralees Nadarajah\\
School of Mathematics\\
University of Manchester\\
Manchester M13 9PL, UK
\end{center}
\vspace{1.5cm}
{\bf Abstract:}~~We give the distribution function of $M_n$,  the maximum of a
sequence of $n$ observations from an autoregressive process of order 2.
Solutions are first given in terms of repeated integrals
and then for the case, where the underlying random variables
are absolutely continuous.
When the correlations are positive,
\begin{eqnarray*}
P(M_n  \leq x)\ =a_{n,x},
\end{eqnarray*}
where
\begin{eqnarray*}
a_{n,x}= \sum_{j=1}^\infty \beta_{jx} \ \nu_{jx}^{n} = O \left( \nu_{1x}^{n} \right),
\end{eqnarray*}
where $\{\nu_{jx}\}$ are the eigenvalues of a non-symmetric Fredholm kernel,
and $\nu_{1x}$ is the eigenvalue of maximum magnitude.
The weights $\beta_{jx}$ depend on the $j$th left and right eigenfunctions of the kernel.

These results are large deviations expansions for estimates, since
the maximum need not be standardized to have a limit.
In fact such a limit need not exist.

\noindent
{\bf Keywords:}~~Autoregressive process; Fredholm kernel; Maximum.

\section{Introduction and summary}
\setcounter{equation}{0}

Many authors have considered extreme value theory for moving average processes,
see Rootz\'{e}n (1978), Leadbetter {\it et al.} (1983, page 59),
Davis and Resnick (1985), Rootz\'{e}n (1986), O'Brien (1987), Resnick (1987, page 239),
Davis and Resnick (1989), Park (1992),
Hall (2002), Hall (2005) and Kl\"{u}ppelberg and Lindner (2005).
However, the results either give the limiting extreme value distributions or assume that the errors come from a specific class
(e.g. integer-valued, exponential type, heavy tailed, light tailed, etc).
We are aware of no work giving the {\it exact} distribution of the maximum of moving average processes.

This paper applies a powerful new method for giving the {\it exact}
distribution of extremes of $n$ correlated observations as weighted
sums of $n$th powers of associated eigenvalues.
The method was first illustrated for a moving average of order 1
 in Withers and Nadarajah (2009a)
 and an  autoregressive process of order 1  in Withers and Nadarajah (2009b).

Let  $\{e_i\}$ be independent and identically distributed random variables from some distribution function $F$ on $R$.
We consider the autoregressive process of order 2,
\begin{eqnarray}
X_i= e_i+ r_1 X_{i-1}+ r_2 X_{i-2}.
\label{ma1}
\end{eqnarray}
 We restrict ourselves to the case where 
\begin{eqnarray}
r_1>0, \ r_1^2+ r_2>0.\label{rr}
\end{eqnarray}
This includes the most important case, $r_1>0, r_2>0$. (When this condition
does not hold the method can be adapted as done in
 Withers and Nadarajah (2009b).)
In Section 2, we give expressions for the distribution function of the maximum
$$
M_n=\max_{i=1}^n X_i,\ n\geq 1,
$$
 in terms of repeated integrals.
This is obtained via the recurrence relationship
\begin{eqnarray}
G_{n}(y) = {\cal K}  G_{n-2}(y),\ y=(y_0,y_1),
\
n\geq 2,
\label{joint}
\end{eqnarray}
where
\begin{eqnarray}
G_{n}(y) &=& P(M_n\leq x,\ X_n\leq y_0,\ X_{n-1}\leq y_1 ) ,
\label{defG}\\
{\cal K} r(y)  &=& E\ \int r( g_{y}( z_1, e_{1},e_0), dz_1),
\label{calK}\\
g_{y}( z_1, e_{1},e_0)  &=& \min_{j=1,2} g_j,\nonumber\\
 g_1  &=&  (y_{0x}-e_0-r_1 e_1-r_1r_2z_1)/(r_1^2+r_2),\
              g_2= (y_{1x}- e_{1}- r_2 z_1)/r_1,\label{gee} \\
 y_{ix} &= &\min(y_i,x), 
\label{g}
\end{eqnarray}
 $I(A)=1$ or 0 for $A$ true or false
and dependency on $x$ is suppressed except in $y_{ix}$.
So, $ {\cal K}$ is a linear  integral operator depending on $x$.
For (\ref{joint}) to work at $n=2$ we define $M_0=-\infty$ so  that  
\begin{eqnarray}
G_{0}(y)=P(X_0\leq y_0, X_{-1}\leq y_1)=H(y)\mbox{ say}.\label{H}
\end{eqnarray}

In Section 3, we consider the case when $F$ is absolutely continuous
with density $f(x)$ with respect to Lebesque measure.
In this case we show that corresponding to ${\cal K}$ is a Fredholm kernel $K(y,z)$.
We give a solution in terms of its eigenvalues and eigenfunctions.
This leads easily to the asymptotic results stated in the abstract.

Our expansions for $P(M_n\leq x)$ for fixed $x$ are large deviation results.
If $x$ is replaced by $x_n$ such that  $P(M_n\leq x_n)$ tends to the generalized extreme value distribution function,
then the expansion still holds, but not the asymptotic expansion
in terms of a single eigenvalue, since this may approach 1 as $n\rightarrow \infty$.

For $a,b$ functions on $R^2$, set 
$\int a=\int a(y)dy=\int_{R^2} a(y)dy$ and 
similarly for
$\int ab$.

\section{Solutions using repeated integrals}
\setcounter{equation}{0}

\begin{theorem}
 $G_n$ of (\ref{defG}) satisfies the recurrence relation (\ref{joint}) in terms
of the integral operator ${\cal K}$ of (\ref{calK}).
\end{theorem}

\noindent
{\bf Proof:} 
Set  
$$c_{y}( X,e) = \min((y_{0x}-e - r_2 X)/r_1,  y_{1x}).
$$
For $n\geq 1$, $G_n$ of (\ref{defG}) satisfies
\begin{eqnarray*}
G_{n}(y)
&=& P(M_n\leq x,\ X_n\leq y_{0x},\ X_{n-1}\leq y_{1x})\\
&=& P(M_{n-1}\leq x, 
e_n+ r_1 X_{n-1}+ r_2 X_{n-2}\leq y_{0x},\
         X_{n-1}\leq y_{1x}),\ n\geq 1,
\\
&=&
P(M_{n-1}\leq x,  X_{n-1}\leq c_{y}( X_{n-2},e_n)) \mbox{ since }r_1>0
\\
&=&
P(M_{n-2}\leq x,  e_{n-1}+ r_1 X_{n-2}+ r_2 X_{n-3}
\leq c_y( X_{n-2},e_n)),\ n\geq 2,\\
&=&
P(M_{n-2}\leq x,  X_{n-2}\leq g_y( X_{n-3}, e_{n-1},e_n))={\cal K}G_{n-2}(y).
\end{eqnarray*}
So for $n\geq 2$, (\ref{joint}) holds.
This ends the proof.
$\Box$

Our goal is to determine $u_n=P(M_n\leq x)= G_{n}(\pmb{\infty})$
where $\pmb{\infty}=(\infty,\infty)$.
Our next result gives these in terms of
\begin{eqnarray*}
a_n= [{\cal K}^n  H(y)]_{y=\pmb{\infty}},\
b_n= [{\cal K}^n  H(y_{0x},y_1)]_{y=\pmb{\infty}},\
 n\geq 0.
\end{eqnarray*}
For example,
\begin{eqnarray*}
a_0 &=& 1, \ a_1= E\ \int H(g_{\pmb{ \infty}}( s, e_{1},e_0), ds),
\\
b_0 &=& H(x,\infty)=P(X_0\leq x),
\
b_1= E\ \int G_1(g_{\pmb{ \infty}}( s, e_{1},e_0), ds),
\end{eqnarray*}
where
$$g_{\pmb{ \infty}}( s, e_{1},e_0) 
=\min\{ (x-e_0-r_1e_1-r_1r_2s)/(r_1^2+r_2), 
       (x-e_1-r_2s)/r_1\}.
$$
\begin{theorem}
\begin{eqnarray}
u_{2n} = a_n,\ u_{2n-1} = b_n, \ n\geq 0.
\label{pos}
\end{eqnarray}
\end{theorem}

\noindent
{\bf Proof:}
 by Theorem 2.1, for $n\geq 0$
\begin{eqnarray*}
G_{2n}(y)= {\cal K}^n  G_0(y),\
 G_{2n+1}(y)= {\cal K}^n  G_1(y).
\end{eqnarray*}
Also 
$$G_0(y)=H(y),\ G_1(y)=H(y_{0x},y_1).$$
Putting $y=\pmb{\infty}$ gives (\ref{pos}).
$\Box$ \\
Note that
$$G_2(y)=H(y_{0x},y_{1x}).$$
\section{The case of $F$ absolutely continuous}
\setcounter{equation}{0}

Our solution Theorem 2.2 does not tell us how
$u_n$ behaves for large $n$.
Also calculating $a_n$ requires repeated integration.
Here we give another solution that overcomes these problems,
using Fredholm integral theory given in Appendix A of 
Withers and Nadarajah (2009a), referred to below as ``the appendix''.

\begin{theorem}
Suppose that $F$ has first and second derivatives $f$ and $f_{.1}$.
Suppose that 
\begin{eqnarray}
\gamma_y(z)r(z)\rightarrow 0\mbox{ as }z_1\rightarrow \pm \infty
\label{inf}
\end{eqnarray}
\begin{eqnarray*}
\mbox{where }z &=& (z_0,z_1),\\
 \gamma_y(z) &=&  (r_1+r_2/r_1)\gamma_{y1}(z) +r_1\gamma_{y2}(z),\\
 \gamma_{y1}(z) &=& \int_{z_0<\beta_y(w_0)} f(c_{y}(w_0,z))f(w_0)dw_0,\\
\beta_y(w_0) &=&  (y_{0x}-r_1y_{1x}-w_0)/r_2,\\ 
c_y(w_0,z) &=& [y_{0x}-r_1r_2z_1-(r_1^2+r_2)z_0-w_0]/r_1,\\
\gamma_{y2}(z) &=& f(y_{1x}-r_2z_1-r_1z_0)\ F(\delta_y(z_0)),\\
\delta_y(z_0)  &=& y_{0x}-r_1y_{1x}-r_2z_0.
\end{eqnarray*}
Then we can write (\ref{calK}) in the form
\begin{eqnarray}
{\cal K} r(y) = \int K(y,z) r(z)dz
\label{kernel}
\end{eqnarray}
where
\begin{eqnarray*}
K(y,z) &=& (r_1+r_2/r_1)r_2\int^{\delta(z_0)} f(\alpha_1(w_0,y,z))f(w_0)dw_0
+ r_1r_2F(\delta(z_0))f_{.1}(\alpha_2(y,z)),\\
\alpha_1(w_0,y,z) &=& [y_{0x}-(r_1^2+r_2)z_0-r_1r_2z_1-w_0]/r_1,\\
\alpha_2(y,z) &=&  y_{1x}-r_1z_0-r_2z_1.
\end{eqnarray*}
\end{theorem}
PROOF
 Set
$$h_y(e_0,z_1)
= [(r_1^2+r_2) y_{1x} -r_1y_{0x} +r_1e_0-r_2^2z_1]/r_2.
$$
Then for $g_i$ of (\ref{gee}),
$$
g_1\leq g_2 \iff e_1\leq h_y(e_0,z_1).$$
So ${\cal K}r(y)=I_1+I_2$ where
\begin{eqnarray*}
I_1 &=& \int f(w_1)dw_1\int_{w_1\leq h_y(w_0,z_1)} f(w_0)dw_0\int r(g_1,dz_1),\\
I_2 &=& r_1\int f(w_1)dw_1\int_{g_2\geq \beta_y(w_0)} f(w_0)dw_0\int r(g_2,dz_1),\\
\mbox{since }r_1\beta_y(w_0) &=& y_{1x}-r_2z_1- h_y(w_0,z_1).
\end{eqnarray*}
%So $r_2\beta(w_0)=y_{0x}-r_1y_{1x}-w_0$ does not depend on $z_1$.
Also $g_1=w_1\iff w_1=c_y(w_0,z_1,g_1)$. So transforming from $g_1$ to $z_0$,
$$I_1=\int f(w_0)dw_0\int A(w_0,z_1)$$
where 
\begin{eqnarray*}
A(w_0,z_1) &=& \int_{h_y(w_0,z_1)}^\infty dw_1f(w_1) r(z_0,dz_1) 
=(r_1+r_2/r_1)\int^{\beta_y(w_0)} f(c(w_0,z))r(z_0,dz_1).\\
\mbox{So }I_1 &=& (r_1+r_2/r_1)\int dz_0\int r(z_0,dz_1) \gamma_{y1}(z).
\end{eqnarray*}
Also 
$$I_2=\int dw_0f(w_0)\int \int^{h_y(w_0,z_1)}dw_1 r(g_2,dz_1)f(w_1)
$$
where $r_1g_2=y_{1x}-w_1-r_2z_1$, 
that is, $w_1=y_{1x}-r_2z_1-r_1g_2$.
  So transforming from $g_2$ to $z_0$,
$$I_2=r_1\int dw_0f(w_0)\int \int_{\beta(w_0)}^\infty dz_0
f(y_{1x}-r_2z_1-r_1z_0) r(z)
=r_1\int dz_0 \int \gamma_{y2}(z) r(z_0,dz_1).
$$
So ${\cal K}r(y)=I_1+I_2=\int dg\int \gamma_y(z)r(z_0,dz_1)$.
%Switching to $z_0=g, z_1=z$ and 
Integrating by parts, (\ref{inf}) gives
$${\cal K}r(y)=\int \gamma_y(z)r_{.1}(z)dz
=\int K_0(y,z)r(z)dz$$
where $r_{.1}(z)=-(\partial/\partial z_1)r(z)$
and $ K_0(y,z)=-\gamma_{y.1}(z)=-(\partial/\partial z_1)\gamma_y(z)
$
by (\ref{inf}).
Also $K_0(y,z)=K(y,z)$.
This ends the proof.
$\Box$ \\

We now assume that
\begin{eqnarray}
0< \int\int  K(y,z)K(z,y)dydz <\infty. \label{as1}
\end{eqnarray}
So $K(y,z)$ is a (non-symmetric) Fredholm kernel with respect
to Lebesgue measure, allowing
the Fredholm theory of the appendix
to be applied, in particular the
functional forms of the Jordan form and singular value decomposition.

Let  $\{\lambda_{j},r_{j},l_{j}:\ j\geq 1 \}$ be the eigenvalues and associated
right and left eigenfunctions of ${\cal K}$ ordered so that $|\lambda_{j}|\geq |\lambda_{j+1}|$. If $\{\lambda_{j} \}$ are real then 
$\{r_{j},l_{j} \}$ can be taken as real.
By the appendix referred to, these satisfy
\begin{eqnarray*}
{\cal K} r_{j}=\lambda_{j} r_{j},
\
\bar{l}_{j} {\cal K}=\lambda_{j} \bar{l}_{j},
\
\int r_{j}\bar{l}_{k}=\int_{R^2} r_{j}(y)\bar{l}_{k}(y)dy
=\delta_{jk},
\end{eqnarray*}
where $\delta_{jk}$ is the Kronecker function.
So, $\{ r_{j}(y), l_{k}(y)\}$ are biorthogonal functions with respect to Lebesgue measure.

We now assume that 
\begin{eqnarray}
K(y,z)\mbox{ has diagonal Jordan form.} \label{as2}
\end{eqnarray}
(This holds, for example, when the eigenvalues are distinct.
This will generally be the case for our applications.)
The functional equivalent of the Jordan form is, by (3.6) of Withers and Nadarajah (2008b), 
\begin{eqnarray*}
K(y,z) = \sum_{j=1}^\infty \lambda_j r_j(y)\bar{l}_j(z).
\end{eqnarray*}
This implies that
\begin{eqnarray*}
K_n(y,z)={\cal K}^{n-1} K(y,z) = \sum_{j=1}^\infty \lambda_j^n r_j(y)\bar{l}_j(z)
\end{eqnarray*}
where ${\cal K}^n$ is the operator corresponding to the iterated kernel
$ K_n(y,z)$.
By (A.8) of Withers and Nadarajah (2009a) with $\mu$ Lebesgue measure on $R^2$,
 if ${\cal K} G$ is in $L_2(R^2)$
then 
\begin{eqnarray}
{\cal K}^n G(y)=\sum_{j=1}^\infty B_j(G)\ r_j(y)\lambda_j^n,
\
n\geq 1,
\label{n}
\end{eqnarray}
where $B_j(G)=\int_{R^2} G\bar{l}_j$.
Putting $y=\pmb{ \infty}$ and $G=G_0,G_1$ in (\ref{n}) gives

\begin{theorem} %3.2
Suppose that (\ref{as1}) (\ref{as2}) hold. Then
for $B_j$ of (\ref{n}) and $n\geq 1$
\begin{eqnarray*}
a_n=\sum_{j=1}^\infty r_j(\pmb{ \infty}) B_j(H) \lambda_j^n,\
b_n=\sum_{j=1}^\infty r_j(\pmb{ \infty}) B_j(G_1) \lambda_j^n.
\end{eqnarray*}
\end{theorem}

\begin{corollary}
Suppose that the eigenvalue $\lambda_1$ of largest magnitude has multiplicity $M$. Then under the assumptions of Theorem 3.2,
\begin{eqnarray}
a_n  = B(H) \lambda_1^{n}(1+\epsilon_{n}),\
b_n  = B(G_1) \lambda_1^{n}(1+\epsilon_{n}),\ n\geq 1,
\label{deltav}
\end{eqnarray}
where $\epsilon_{n}\rightarrow 0$ exponentially as $n\rightarrow \infty$
and
\begin{eqnarray*}
B(G) =\sum_{j=1}^M  r_j(\pmb{\infty}) B_j(G).
\end{eqnarray*}
So, for $n\geq 1$, by (\ref{pos})
\begin{eqnarray*}
u_{2n} =
  B(H) \lambda_1^{n}(1+\epsilon_{n}),\
u_{2n+1} =  B(G_1) \lambda_1^{n}(1+\epsilon_{n}).
\end{eqnarray*}
\end{corollary}

\subsection{A numerical solution}

We now give a numerical method for
obtaining the eigenvalues and
 $r_j(\pmb{ \infty})$ and $B_j(G)$ needed for Theorem 3.2.

Consider the $j$th eigenvalue and eigenfunctions
 $\lambda=\lambda_j,\ r= r_{j}, \ l=l_j$.
The  right and left eigenfunctions satisfy 
$$\lambda r={\cal K}r,\ \lambda \bar{l}= \bar{l}{\cal K},$$
 that is
\begin{eqnarray}
\lambda r(y)= {\cal K} r(y) = \int_{R^2} K(y,z)r(z)dz,\
\lambda \bar{l}(z)=  \bar{l}(z){\cal K}  = \int_{R^2}  \bar{l}(y)K(y,z)dy.
\label{kay}
\end{eqnarray}
Let us approximate an integral over ${R^2}$ by
by Gaussian quadrature 
(see for example Section 25.4 of  Abramowitz and Stegun (1964)), say
\begin{eqnarray*}
\int_{R^2} a(z)dz \approx \sum_{j=1}^r w_j a(z_j),
\end{eqnarray*}
where $\{z_1,\cdots, z_r\}$ are given points in  ${R^2}$ and
$\{w_1,\cdots, w_r\}$ are given weights.
Let $ {\bf R, L}$ denote the $r$-vectors with $k$th elements $r(z_k),l(z_k)$.
Let $ {\bf K}$ denote the $r\times r $ matrix with
 $(i,j)$th element $K(z_i,z_j)$.
Then we can write (\ref{kay}) as 
$$
\lambda {\bf R}\approx 
{\bf K} {\bf R},\ \lambda  \bar{\bf L}\approx  \bar{\bf L}{\bf K}.
$$
So to this order of approximation,
the eigenvalues are  just those of  $ {\bf K}$,
and  $r(z_k),l(z_k)$ are just the  $k$th elements of
the right and left eigenvectors of  $ {\bf K}$ corresponding to $\lambda$.

Also
$$r(\pmb{ \infty})=\lambda^{-1}\int K(\pmb{\infty},z)r(z)dz
\approx \lambda^{-1}\sum_{k=1}^r w_k K(\pmb{\infty},z_k) r(z_k)$$
and 
$$B_j(G)=\int_{R^2} G\bar{l}_j\approx \sum_{k=1}^r w_k G(z_k)\bar{l}_j(z_k).$$


\begin{thebibliography}{999}


\bibitem{} 
Abramowitz, M. and Stegun, I. A. (1964) {\it Handbook of mathematical functions.} U.S. Department of Commerce, National Bureau of Standards, Applied Mathematics Series  {\bf 55}.
% Abramowitz and Stegun (1964)

\bibitem{}
Davis, R. A. and Resnick, S. I. (1985).
Limit theory for moving averages of random variables with regularly varying tail probabilities.
{\it Annals of Probability}, {\bf 13}, 179--195.


\bibitem{}
Davis, R. A. and Resnick, S. I. (1989).
Basic properties and prediction of max--ARMA processes.
{\it Advances in Applied Probability}, {\bf 21}, 781--803.

%\bibitem{} Dieudonné, Jean (1953). On biorthogonal systems. {\it Michigan Math. J.}, {\bf 2}, no. 1, 7--20 [1].

\bibitem{}
Hall, A. (2002).
Extremes of integer--valued moving average models with regularly varying tails.
{\it Extremes}, {\bf 4}, 219--239.


\bibitem{}
Hall, A. (2005).
Extremes of integer--valued moving average models with exponential type tails.
{\it Extremes}, {\bf 6}, 361--379.



\bibitem{}
Kl\"{u}ppelberg, C. and Lindner, A. (2005).
Extreme value theory for moving average processes with light-tailed innovations.
{\it Bernoulli}, {\bf 11}, 381--410.


\bibitem{}
Leadbetter, M. R., Lindgren, G. and Rootz\'{e}n, H. (1983).
{\it Extremes and Related Properties of Random Sequences}.
Springer-Verlag, New York.


\bibitem{}
O'Brien, G. L. (1987).
Extreme values for stationary and Markov sequences.
{\it Annals of Probability}, {\bf 15}, 281--291.


\bibitem{}
Park, Y. S. (1992).
Extreme value of moving average processes with negative binomial noise distribution.
{\it Journal of the Korean Statistical Society}, {\bf 21}, 167--177.


\bibitem{}
Resnick, S. I. (1987).
{\it Extreme Values, Regular Variation, and Point Processes}.
Springer-Verlag, New York.


\bibitem{}
Rootz\'{e}n, H. (1978).
Extremes of moving averages of stable processes.
{\it Annals of Probability}, {\bf 6}, 847--869.


\bibitem{}
Rootz\'{e}n, H. (1986).
Extreme value theory for moving average processes.
{\it Annals of Probability}, {\bf 14}, 612--652.



\bibitem{}
Withers, C. S. (1975).
Fredholm theory for arbitrary measure spaces.
{\it Bulletin of the Australian Mathematical Society}, {\bf 12}, 283--292.


\bibitem{}
Withers, C. S. (1978).
Fredholm equations have uniformly convergent solutions.
{\it Journal of Mathematical Analysis and Its Applications}, {\bf 64}, 602--609.



\bibitem{}
Withers, C. S. and Nadarajah, S. (2008a).
The $n$th power of a matrix and approximations for large $n$.
{\it New Zealand Journal of Mathematics}, {\bf 38}, 171--178. %138cor/138corb

\bibitem{}
Withers, C. S. and Nadarajah, S. (2008b).
Fredholm equations for non-symmetric kernels with applications to iterated integral operators.
{\it Applied Mathematics and Computation}, {\bf 204}, 499--507. %138cor/138cora.tex

\bibitem{}
Withers, C. S. and Nadarajah, S. (2009a).
The distribution of the maximum of a first order moving average: the continuous case.
{\it Technical Report}, Applied Mathematics Group, Industrial Research Ltd., Lower Hutt, New Zealand.
Available on-line at {\sf http://arxiv.org/abs/0802.0523}.
\bibitem{}

Withers, C. S. and Nadarajah, S. (2009b).
The distribution of the maximum of a first order autoregressive process: the continuous case.
{\it Technical Report}, Applied Mathematics Group, Industrial Research Ltd., Lower Hutt, New Zealand.

\end{thebibliography}
\end{document}